\documentstyle[leqno]{article}
   
\font\tengothic=eufm10
\newfam\Gothic
\textfont\Gothic=\tengothic
\def\go{\fam\Gothic\tengothic}

  \newtheorem{theorem}{Theorem}
  \newtheorem{lemma}{Lemma}
  \newtheorem{corollary}{Corollary}
  
\begin{document}
\author{D. A. Shmel'kin.}
\title{
First Fundamental Theorem of Invariant Theory
for covariants of classical groups.}
\maketitle

{\bf Abstract.} Let $U(G)$ be a maximal unipotent subgroup
of one of classical groups $G=GL(V),O(V),Sp(V)$.
Let $W$ be a direct sum of copies of $V$ and its dual $V^*$.
For the natural action $U(G):W$, we describe a minimal system
of homogeneous generators for the algebra of $U(G)$-invariant
regular functions on $W$. For $G=GL(V)$, we also describe
the syzygies among these generators in some particular cases.

\vspace{0.5cm}

\section{Main theorem.}\label{main}

   Let $V$ be a finite-dimensional vector space  over an
   algebraically closed field ${\bf k}$ of characteristic zero.
   Let $H\subseteq GL(V)$ be an algebraic subgroup. For
   any $l\in{\bf N}$, we denote by
   $lV$ the direct sum of $l$ copies of $V$;
   similarly, we define $mV^*$ for any $m\in {\bf N}$.
   Consider the natural action of $H$ on $W=lV\oplus mV^*$ and assume that
   the algebra ${\bf k}[W]^H$ of invariants is finitely generated
   for any $l,m$.  Then First Fundamental Theorem of Invariant
   Theory of $H$ refers to a description of a minimal system
   of homogeneous generators of ${\bf k}[W]^H$ for all $l,m$.

   Such a description exists when $H$ is {\it classical}, i.e.,
   $H$ is one of groups
   $GL(V)$, $SL(V)$, $O(V)$, $SO(V)$, $Sp(V)$ (see e.g.
   ~\cite[$\S$9]{vipo}).

   Let now $G$ be one of the groups $GL(V),O(V),Sp(V)$;
   let $U(G)$ be a maximal unipotent subgroup of $G$.
   By ~\cite[Theorem~3.13]{vipo}, the algebra 
   ${\bf k}[W]^{U(G)}$
   is finitely generated.
   Also the invariants of $U(G)$ are linear
   combinations of highest vectors of 
   irreducible factors for $G$-module ${\bf k}[W]$.
   So the $U(G)$-invariants are the $G$-covariants
   and the First Fundamental Theorem for  
   covariants of $G$ means that for the invariants
   of $H=U(G)$.
   
   Using a result (and some ideas) of
   ~\cite{ho},
   we prove in this paper First Fundamental
   Theorem for covariants of each of the above
   classical $G$.
   
   Note that for $G=Sp(V),O(V)$, $V$ and $V^*$ are 
   isomorphic as $G$-modules,
   hence, as $U(G)$-modules. Therefore we may assume 
   $m=0$ in these cases.
   Elements of
$\wedge^kV^*\subseteq\otimes^kV^*\subseteq{\bf k}[kV]$, $k\leq \dim V$
are said to be {\it multilinear anisymmetric} functions as well as
their analogs in ${\bf k}[kV^*]$.

\begin{theorem}\label{generators}
The algebra ${\bf k}[W]^{U(G)}$ is generated by
the subalgebra ${\bf k}[W]^G$ and multilinear antisymmetric invariants.
Moreover, a set ${\cal M}$ described below is a
minimal system of homogeneous generators of ${\bf k}[W]^{U(G)}$.
\end{theorem}

We describe a minimal system ${\cal M}$ of homogeneous
generators of ${\bf k}[W]^{U(G)}$ in a coordinate form.
Set $n=\dim V$, choose a basis of $V$ and
denote by $\overline{V}$ the corresponding $n\times l$-matrix of coordinates
on $lV$. Similarly, denote by $\overline{V^*}$ the $m\times n$-matrix
of coordinates on $mV^*$, in the dual basis of $V^*$. A minor
of order $k$ of a matrix is said to be {\it left},
if it involves the first $k$ columns.
Analogeously, we call it {\it lower}, if it involves the
last $k$ rows. 

\vspace{0.2cm}

A) Let $G=GL(V)$ and define $U(GL)=U(GL(V))$
to be the subgroup of the strictly upper triangular matrices,
in the above basis. Then ${\cal M}$ is:

\noindent$\bullet$ the matrix elements of the product $\overline{V^*}\overline{V}$

\noindent$\bullet$ the lower minor determinants of order $k$
of $\overline{V}$, $k=1,\cdots,{\rm min}\{l,n\}$.

\noindent$\bullet$ the left minor determinants of order $p$
of $\overline{V^*}$, $p=1,\cdots,{\rm min}\{m,n\}$.

Let $T(GL)$ be the diagonal matrices in the above basis.
Then $T(GL)$ is a maximal torus of $G$ normalizing $U(GL)$.
The pair $T(GL),U(GL)$ defines a system of simple roots of $T(GL)$.
Here and in what follows, we use the enumeration of simple roots
of simple groups as in ~\cite{ov} and
denote by $\varphi_1,\cdots,\varphi_n$ the fundamental weights.
The torus $T(GL)$ acts on ${\bf k}[W]^{U(GL)}$ and  the 
elements of ${\cal M}$ are weight vectors of $T(GL)$. The set of their
degrees and weights is (for $l,m\geq n$):

$$(2,0),
(1,\varphi_1),(2,\varphi_2),\cdots,(n,\varphi_n),$$
$$(1,\varphi_{n-1}-\varphi_n),(2,\varphi_{n-2}-\varphi_n),\cdots,
(n-1,\varphi_1-\varphi_n),(n,-\varphi_n).$$

\vspace{0.2cm}

Furthermore,
let $Q$ be a bilinear
symmetric (antisymmetric) form having in the above basis a matrix
with $\pm1$ on the secondary diagonal and with zero entries
outside it. Define $G=O(V)$ ($G=Sp(V)$) to be the stabilizer of
this form. Then $U(G)=G\cap U(GL)$ is a maximal
unipotent subgroup in $G$.
Moreover, set $T(O)=T(GL)\cap SO(V)$, $T(Sp)=T(GL)\cap Sp(V)$.
Then $T(G)$ is a maximal torus of $G$ of rank $r=[${\large
$\frac{n}{2}$}$]$.
Denote by $\varphi_1,\cdots,\varphi_r$ the fundamental weights
of $T(G)$ with respect to $U(G)$.
For $x\in W=lV$, denote by $v_i$ the projections of
$x$ on the $i$-th $V$-factor, $i=1,\cdots,l$.

\vspace{0.2cm}

B) Let $n=2r+1,G=O(V)$. Then  ${\cal M}$ is:

\noindent$\bullet$ $Q(v_i,v_j),1\leq i\leq j\leq l$,

\noindent$\bullet$ the lower minor determinants of order $k$
of $\overline{V}$, $k=1,\cdots,{\rm min}\{l,n\}$,

The set of degrees and weights of the above generators is
(for $l\geq n$):

$$(2,0),(1,\varphi_1),\cdots,(r-1,\varphi_{r-1}),
(r,2\varphi_r),(r+1,2\varphi_r),\cdots,(n-1,\varphi_1),(n,0).$$

\vspace{0.2cm}

C) Let $G=Sp(V)$.
Then ${\cal M}$ is:

\noindent$\bullet$ $Q(v_i,v_j),1\leq i<j\leq l$,

\noindent$\bullet$ the lower minor determinants of order $k$
of $\overline{V}$, $k=1,\cdots,{\rm min}\{l,r\}$.

The set of degrees and weights of the above generators
is (for $l\geq r$):

$$(2,0),(1,\varphi_1),(2,\varphi_2),\cdots,
(r,\varphi_r).$$

Note that the lower minor determinants of order $k$ of
$\overline{V}$ with $k>r$ are
$U(Sp)$-invariant, too. It is not hard to check that these
can be expressed in the above generators.

\vspace{0.2cm}

D) Let $n=2r,G=O(V)$.
Then ${\cal M}$ is:

\noindent$\bullet$ $Q(v_i,v_j),1\leq i\leq j\leq l$,

\noindent$\bullet$ the lower minor determinants of order $k$
of $\overline{V}$, $k=1,\cdots,{\rm min}\{l,n\}$,

\noindent$\bullet$ for $l\geq r$, the minor determinants of order $r$,
involving the $r$-th row and the last $r-1$ rows of $\overline{V}$.

The set of degrees and weights of the above generators is
(for $l\geq n$):

$$ (2,0),(1,\varphi_1),\cdots,(r-2,\varphi_{r-2}),
(r-1,\varphi_{r-1}+\varphi_r),(r,2\varphi_{r-1}),(r,2\varphi_r),$$
$$(r+1,\varphi_{r-1}+\varphi_r),\cdots,(n-1,\varphi_1),(n,0).$$

\section{Proof of Theorem 1.}\label{pf}

 First we state a result of ~\cite{ho} that
 is a starting point of our proof.
 We keep the notation of loc.cit.
 but consider a slightly more general setting.
 
 Let $W$ be a finite dimensional ${\bf k}$-vector
 space. Denote by
 $${\go gr}={\go gr}_{(2,0)}\oplus{\go gr}_{(1,1)}
 \oplus{\go gr}_{(0,2)}\subseteq {\rm End}{\bf k}[W]$$
 the linear subspace of differential operators
 with the prescribed by the index degree and order.
 Namely, ${\go gr}_{(2,0)}$ are the homogeneous
 regular functions on $W$ of degree 2 acting on
 ${\bf k}[W]$ by multiplication; ${\go gr}_{(0,2)}$
 are the constant coefficients differential
 operators of order 2; ${\go gr}_{(1,1)}$ is
 nothing but the Lie algebra ${\go gl}(W)$.
 
 Clearly, ${\go gr}$ is a Lie subalgebra in
 ${\rm End}{\bf k}[W]$, and moreover,
 ${\go gr}$ is isomorphic to 
 ${\go sp}(W\oplus W^*)$, with respect to
 the natural symplectic form on $W\oplus W^*$.

 Assume now that $G\subseteq GL(W)$ is a reductive subgroup.
 Then $G$ acts on ${\go gr}$; consider the invariants:
 $$
 \Gamma'={\go gr}^G, 
 {\Gamma'}_{(2,0)}={\go gr}_{(2,0)}^G,
 {\Gamma'}_{(1,1)}={\go gr}_{(1,1)}^G, 
 {\Gamma'}_{(0,2)}={\go gr}_{(0,2)}^G.$$
 Clearly, 
 $\Gamma'= {\Gamma'}_{(2,0)}\oplus
 {\Gamma'}_{(1,1)}\oplus {\Gamma'}_{(0,2)}$
 is also a Lie subalgebra in ${\rm End}{\bf k}[W]$.
 
  Let ${\bf k}[W]=\bigoplus\limits_{k=1}^{\infty} I_k$ be 
  the decomposition of $G$-module ${\bf k}[W]$ into
  isotypic components. Let $I$ be one of $I_k$.
  Clearly, $I$ is stable under the action of
  $\Gamma'$.

 \begin{theorem}\label{irreduce}
 {\rm (~\cite[Theorem~8]{ho})}
 Assume that the algebra ${\bf k}[W\oplus W^*]^G$
 of invariants is generated by elements of degree 2.
 Then $I$ is an irreducible joint $(G,\Gamma')$-module.
 \end{theorem} 
    
 By the First Fundamental Theorem for the classical groups,
 the assumption of Theorem ~\ref{irreduce} holds
 for the pairs $(G,W)$ from section ~\ref{main}.
 For these particular cases the above Theorem
 is (a part of) Theorem 8 of ~\cite{ho}. However,
 one can see that the proof in loc.cit. works 
 whenever the assumption of Theorem ~\ref{irreduce} holds.  
 
 Note that for the classical $(G,W)$ we have:
 $\Gamma'_{(1,1)}={\go gl}_l\oplus{\go gl}_m$, 
  $$\Gamma'\cong {\go gl}_{l+m},\;{\rm if}\;G=GL(V),
    \Gamma'\cong {\go sp}_{2l},\;{\rm if}\;G=O(V),
    \Gamma'\cong {\go o}_{2l},\;{\rm if}\;G=Sp(V).$$

  We now show that Theorem~\ref{irreduce}
  reduces Theorem ~\ref{generators} to
  a more simple statement.
  The below reasoning is an analog of that
  from the proof of Theorem 9 in loc.cit.
  
  Clearly, $I$ is a homogeneous submodule
  of ${\bf k}[W]$; denote by $I^{min}$ the subspace
  of the elements of $I$ of minimal degree.
  Let $A\subseteq {\bf k}[W]^U$ be the subalgebra
  generated by ${\cal M}$.
  Let $Z\subseteq {\bf k}[W]$ be the
  $G$-submodule generated by $A$. 
  Then Theorem ~\ref{generators} can be
  reformulated as follows: $Z={\bf k}[W]$.
  Assume that $X=Z\cap I^{min}$ is nonzero.
  
  Since the system ${\cal M}$ of generators
  of $A$ is symmetric with respect to permutations
  of isomorphic $G$-factors of $W$, $A$ is
  $GL_l\times GL_m$-stable, i.e., 
  $\Gamma'_{(1,1)}$-stable.
  Hence, $Z$ and $X$ are 
  stable with respect to both $G$ and
  $\Gamma'_{(1,1)}$.
  
  Let $R,R_{(2,0)}$ etc. be the subalgebras in
  ${\rm End}{\bf k}[W]$  generated by
  $\Gamma',{\Gamma'}_{(2,0)}$ etc.
  Consider $R$ as a representation of the universal
  enveloping algebra of $\Gamma'$.
  Using the PBW theorem, we obtain  
  \begin{equation}\label{pbw}
  R=R_{(2,0)}R_{(1,1)}R_{(0,2)}.
  \end{equation}
  Differentiating a polynomial, we decrease 
  its degree; hence, ${\Gamma'}_{(0,2)}I^{min}=0$.
  Therefore $R_{(0,2)}X=X$. Moreover, 
  since $X$ is $\Gamma'_{(1,1)}$-stable, we have by (~\ref{pbw}):
  $RX=R_{(2,0)}X={\bf k}[W]^GX$. 
  On the other hand, $RX$ is a non-zero 
  joint $(G,\Gamma')$-submodule of $I$. 
  By Theorem ~\ref{irreduce}, 
  $I=RX={\bf k}[W]^GX\subseteq Z$.
 
  Thus to prove Theorem ~\ref{generators}, 
  we need to check 
  for any isotypic component $I$: 
  
  \begin{equation}\label{all}
  A\cap I^{min}\neq\{0\}.
  \end{equation}

  Note that it is sufficient to prove 
  Theorem~\ref{generators}
  with $l,m\geq n$, in the case $G=GL(V)$, and
  with $l\geq n, m=0$, in the case $G=O(V),Sp(V)$. 

  Denote by $G^0$ the connected component of the  
  unity of $G$;
  $GL(V)$ and $Sp(V)$ are connected, but for 
  $G=O(V)$, $G^0=SO(V)$.
  Recall that the irreducible finite dimensional 
  $G^0$-modules are
  in one-to-one correspondence
  with their highest weights with respect to $U(G)$ and $T(G)$.
  Denote by $P$ the set of highest weights of irreducible factors
  for $G^0$-module ${\bf k}[W]$. For any graded algebra $B$
  and $t\in {\bf N}$, we denote by $B_t$ the subspace of the elements
  of degree $t$. For any $\chi\in P$ we set:

  $R(\chi)$ is the irreducible
  representation of $G^0$ with highest weight $\chi$

  $I_{\chi}$ is the $R(\chi)$-isotypic component
  of $G^0$-module ${\bf k}[W]$

  $m(\chi)={\rm min}\{t\vert {\bf k}[W]_t\cap
  I_{\chi}\neq 0\}$.

  $n(\chi)={\rm min}\{t\vert A_t\cap I_{\chi}\neq 0\}$.

 By definition, $n(\chi)\geq m(\chi)$.
 For $G=GL(V),Sp(V)$ the condition (~\ref{all})
 is equivalent to $n(\chi)=m(\chi)$ for any  
 $\chi\in P$.

 \begin{lemma}\label{linear}
  For any $\chi\in P,c\in{\bf N}$ we have: $n(c\chi)=cn(\chi)$.
  \end{lemma}

  Denote by ${\go t}$ the Lie algebra of $T(G)$.
  Let ${\cal C}\subseteq {\go t}^*$ be the Weyl chamber corresponding
  to $U(G)$. Consider the set
  $$\Delta=\{\frac{\chi^*}{t}\vert I(\chi)\cap {\bf k}[W]_t\neq 0\}
  \subseteq {\cal C},$$
  where $\chi^*$ denotes the highest weight of the $G^0$-module dual
  to that with highest weight $\chi$.
  By ~\cite{br1}, if ${\bf k}$ is the field ${\bf C}$ of complex numbers,
  then $\Delta$ is the set of rational points in
  the momentum polytope
  for the action of the maximal compact subgroup $K\subseteq G^0$
  on the projective space ${\bf P}(W)$. Further, we set:
  $$\widetilde{\Delta}=
  \{\frac{\chi^*}{t}\vert I(\chi)\cap Z_t\neq 0\} \subseteq \Delta.$$
  Let now $\Phi\subseteq{\go t}^*$ be the
  convex hull over the rational numbers of the weights for the action
  $T(G):W$.

  \begin{lemma} $\widetilde{\Delta}\supseteq\Phi\cap{\cal C}$.\label{delta}
  \end{lemma}

  By definition, we have: $\Delta\subseteq\Phi\cap {\cal C}$.
  Therefore $\Delta=\widetilde{\Delta}=\Phi\cap{\cal C}$
  \footnote{
  For ${\bf k=C}$, one can directly prove for the moment polytope
  $\Delta\otimes{\bf R}=(\Phi\otimes{\bf R})\cap{\cal C}$.}.

  Suppose that ${\bf k}[W]^{U(G)}$
  contains an element of degree $t$ and weight
  $\chi$. Then by definition,
  {\large $\frac{\chi^*}{t}$}$\in \Delta$.
  Hence, the equality $\Delta=\widetilde{\Delta}$
  implies that for some $c\in{\bf N}$ there exists
  an element of $A$ of degree
  $ct$ and weight $c\chi$. Thus $ct\geq n(c\chi)=cn(\chi)$ and
  $t\geq n(\chi)$. In other words, $m(\chi)\geq n(\chi)$, hence
  $m(\chi)=n(\chi)$.
  This completes (modulo Lemmas ~\ref{linear} and ~\ref{delta})
  the proof of Theorem for $G=GL(V),Sp(V)$.

  \vspace{0.2cm}

  Let  $G$ be $O(V)$; to prove Theorem, we apply induction on $n=\dim V$.

  For $n=2$, $U(O)$ is trivial and one can see $A={\bf k}[W]$.

  For $n=3$, $(SO_3,{\bf k}^3)\cong (SL_2,S^2{\bf k}^2)$.
  Since the stabilizer of a point on the dense orbit for the action
  $SL_2:{\bf k}^2$ is a maximal unipotent subgroup in $SL_2$,
  we obtain an isomorphism:
  $${\bf k}[{\bf k}^2+l{\bf k}^3]^{SL_2}\cong {\bf k}[W]^{U(O)}.$$

  \begin{lemma}\label{sl2}
  There exists an isomorphism ${\bf k}[{\bf k}^2+l{\bf k}^3]^{SL_2}\cong
  {\bf k}[(l+1){\bf k}^3]^{SO_3}/(d)$, where $d=Q(v_{l+1},v_{l+1})$.
  \end{lemma}

  {\it Proof:} Consider the morphism
  $$\varphi:{\bf k}^2+l{\bf k}^3\to
  (l+1){\bf k}^3, \varphi(e, Q_1,\cdots,Q_l)=(Q_1,\cdots,Q_l,e^2).$$
  Clearly, $\varphi$ is $SL_2$-equivariant; moreover, $\varphi$ is
  the quotient map with respect to the center of $SL_2$.
  Furthermore, the image of $\varphi$
  is the zero level of $d$. This completes the proof.$\Box$

  Using Lemma~\ref{sl2} and the well-known description of
  ${\bf k}[(l+1){\bf k}^3]^{SO_3}$, one easily deduces the Theorem for $n=3$.

 The step of induction. Assume that Theorem is proven for
 $n-2$. We apply now the Theorem of local structure of
 Brion-Luna-Vust (~\cite{blv}) to get a local version of the assertion
 of Theorem.

 Denote by $x_i^j=\overline{V}_i^j$
 the $i$-th coordinate of $v_j$.
 Set $f=x_n^1\in{\bf k}[W]^U$,  $W_f=\{x\in W\vert f(x)\neq 0\}$.
 Define a mapping:
 $$\psi_f:W_f\to{\go o}(V)^*,\psi(x)(\xi)=\frac{(\xi f)(x)}{f(x)}.$$
 Denote by $P_f$ the stabilizer in $SO(V)$ of the line $\langle f\rangle$.
 Clearly, $P_f$ is a parabolic subgroup in $SO(V)$ containing $U(O)$ and
 $\psi_f$ is $P_f$-equivariant.

 Furthermore, we denote by $e_i^j$ the $i$-th element of the above basis
 in the $j$-th copy of $V$, $x=e_n^1$, $\Sigma=\psi_f^{-1}(\psi_f(x))$.
 Denote by $L$ the stabilizer of $\psi_f(x)$ in $P_f$.
 By ~\cite{blv}, $L$ is a Levi subgroup of $P_f$ and the natural
 morphism
 $$P_f*_L\Sigma\to W_f, (p,\sigma)\to p\sigma$$
 is a $P_f$-equivariant isomorphism. Therefore we have:
 $${\bf k}[W]_f^{U(O)}\cong {\bf k}[W_f]^{U(O)}
 \cong{\bf k}[P_f*_L\Sigma]^{U(O)}.$$
 Also, $P_f=U(O)L$. Hence,
 $${\bf k}[P_f*_L\Sigma]^{U(O)}\cong {\bf k}[\Sigma]^{U(O)\cap L}=
 {\bf k}[\Sigma]^{U(L)},$$
 where $U(L)$ is a maximal unipotent subgroup in $L$. Calculating, we
 have:
 $$(L,\Sigma)\cong(SO_2\times SO_{n-2},\langle e_1^1,e_n^1\rangle_f\times
 (l-1)V).$$
 In other words,
 $${\bf k}[W]^{U(O)}_{x_n^1}\cong{\bf k}[x^1_1,x_n^1,x_1^2,x_n^2,\cdots
,x_1^l,x_n^l]_{x_n^1}\otimes{\bf k}[(l-1){\bf k}^{n-2}]^{U(O_{n-2})}.$$
 The induction hypothesis yields the generators of
 ${\bf k}[(l-1){\bf k}^{n-2}]^{U(O_{n-2})}$. Restricting the elements
 of ${\cal M}$ to $\Sigma$, one can easily deduce:
 \begin{equation}\label{local}
 {\bf k}[W]^{U(O)}_{x_n^1}=A_{x_n^1}.
 \end{equation}

 We return to our consideration of the isotypic components
 of $O(V):{\bf k}[W]$.
  Consider an irreducible representation $\rho$
  of $O(V)$ and its restriction $\rho'$ to $SO(V)$. Here two cases
  occur:

  \noindent$\bullet$ either $\rho'$ is also irreducible,
  $\rho'=R(\chi)$ for some $\chi\in P$

  \noindent$\bullet$ or else $n=2r$,
  $\rho'=R(\chi)+R({\tau}(\chi))$,
  where $\tau$ is the automorphism of the system of simple
  roots of $O(V)$ interchanging the $r-1$-th and the $r$-th roots.

  The latter case is more simple: elements of minimal degree
  in the $\rho$-isotypic component
  are the elements of minimal degree in both $I(\chi)$ and
  $I(\tau(\chi))$ (clearly, $n(\chi)=n(\tau(\chi))$ and
  $m(\chi)=m(\tau(\chi))$). Hence, the above equality $n(\chi)=m(\chi)$
  implies the assertion for such an isotypic component.

  Now consider the former case. Here for any $\rho'=R(\chi)$
  there exist two possibilities for $\rho$: $R(\chi_+)$
  and $R(\chi_-)=R(\chi_+)\otimes{\det}$, where $\det$ is the
  unique nontrivial character of $O(V)$. 
  Moreover, we define explicitly $R(\chi_+)$ and $R(\chi_-)$
  as follows.
  Let $\theta\in
  O(V)\setminus SO(V)$ be an element normalizing $T(O)$ as follows.
  For $n$ odd, $\theta=-Id$. For $n$ even, $\theta$ is the operator
  interchanging the $r$-th and the $r+1$-th elements of the above basis
  and acting trivially on the other basis elements. Note that in both cases
  $\theta(\chi)=\chi$ for any $\chi$, if $n$ is odd and for all $\chi$
  such that $\tau(\chi)=\chi$, if $n$ is even. Now we define
  $R(\chi_{\pm})$ by the condition:
  $$R(\chi_{\pm})(\theta)(u_{\chi})=\pm u_{\chi}$$ for the highest vector
  $u_{\chi}$ of $T(O)$ and $U(O)$ in $R(\chi)$. For instance,
  if $n$ is even, $k\leq r-2$,
  minor determinants of order $k$ of $\overline{V}$ generate
  $R(\varphi_{k+})$  and minor determinants of order $n-k$ generate
  $R(\varphi_{k-})$. Moreover, multiplying two highest vectors of
  ${\bf k}[W]$, we add their weights and multiply as usual
  their $\pm$ subscripts. Thus we control the structure of the
  $O(V)$-module $Z$.

  Define $m(\chi_{\pm}),n(\chi_{\pm})$ as above. 
  Then the condition (~\ref{all}) is equivalent to
  the equality
  $m(\chi_{\pm})=n(\chi_{\pm})$ for any $\chi\in P$
  ($\tau$-invariant for $n$ even).
  For any $\chi=\sum_{i=1}^q k_i\varphi_i, k_q>0$, set
  $t=r-1$, if $q=r, n=2r$ and $t=q$ otherwise.
  Then we have:
  \begin{equation}\label{pm}
  {\rm min}\{n(\chi_+),n(\chi_-)\}=n(\chi),
  \vert n(\chi_+)- n(\chi_-)\vert=n-2t.
  \end{equation}
  Let $g$ be a highest vector of ${\bf k}[W]$
  generating $R(\chi_-)$. Then by (~\ref{local}),
  for some even $j$ we have:$(x_n^1)^jg\in A$.
  Since $(x_n^1)^jg$ generates $R((\chi+j\varphi_1)_-)$,
  we have: $\deg g +j\geq n((\chi+j\varphi_1)_-)$.
  Clearly, $n(\chi+j\varphi_1)=n(\chi)+j$ (see formulae
  (~\ref{B}),(~\ref{D}) below). Hence, (~\ref{pm})
  yields $n((\chi+j\varphi_1)_-)=n(\chi_-)+j$.
  Thus we have $\deg g\geq n(\chi_-)$ implying $m(\chi_-)=n(\chi_-)$.
  The same is true for $\chi_+$. This completes the proof
  of Theorem for $G=O(V)$.$\Box$

  Thus we reduced Theorem ~\ref{generators} to 
  Lemma ~\ref{linear} and Lemma ~\ref{delta}.
  Both are properties of degrees and weights of the given generators,
  and we consider case by case.

\section{Proof of Lemmas ~\ref{linear} and  ~\ref{delta}.}

  {\it Proof of Lemma ~\ref{linear}.}
  
  Recall that $n(\chi)$ is the minimum of degree of the monomials
  in the elements of ${\cal M}$ having weight $\chi$. Clearly, we 
  should not involve the $G$-invariants in a monomial of minimal degree.
  Then for $G=Sp(V),O(V)$ we have no much choice for such
  a monomial and we can write down formulae for $n(\chi)$
  as follows.  Let $\chi=k_1\varphi_1+\cdots+k_r\varphi_r$.

  For $G=Sp(V)$, we have: $n(\chi)=k_1+2k_2+\cdots+rk_r$.

  For $G=O(V), n=2r+1$, $k_r$ is even for $\chi\in P$, and we have:
  \begin{equation}\label{B}
  n(\chi)=k_1+2k_2+\cdots+(r-1)k_{r-1}+r\frac{k_r}{2}.
  \end{equation}
  For $G=O(V), n=2r$, $k_{r-1}+k_r$ is even for $\chi\in P$, and we
  have:
  \begin{equation}\label{D}
  n(\chi)=k_1+2k_2+\cdots+(r-2)k_{r-2}+r\frac{k_{r-1}+k_r}{2}
  -{\rm min}(k_r,k_{r-1}).
  \end{equation}
  These formulae yield the assertion of Lemma.

  Consider the case $G=GL(V)$. The elements of ${\cal M}$
  with non-zero weights have the following weights endowed with degrees:
  $$\alpha_i=\varphi_i,\deg\alpha_i=i, i=1,\cdots,n,$$
  $$\beta_j=\varphi_j-\varphi_n,\deg\beta_j=n-j, j=1,\cdots,n-1,
  \beta_n=-\varphi_n,\deg\beta_n=n.$$
  
  For $\chi=k_1\varphi_1+\cdots+k_n\varphi_n$ consider the presentations
  of $\chi$ as linear combinations of the above weights with positive
  integer coefficients. Define the degree of such a combination as the
  sum of degrees of the summands. We claim that there is a unique presentation
  of minimal degree.

  For any
  $j=1,\cdots,n-1$, all the presentations of $\chi$ contain $k_j$
  summands $\alpha_j$ or $\beta_j$.
  Set $r=[${\large$\frac{n}{2}$}$]$.
  The linear  combination  
  $$\chi'=k_1\alpha_1+\cdots+k_r\alpha_r+
  k_{r+1}\beta_{r+1}+\cdots+k_{n-1}\beta_{n-1}$$ 
  has the minimal
  degree among the linear combinations equal to $\chi$ modulo
  $\langle\varphi_n\rangle$. If $\chi'=\chi$, then this 
  presentation of $\chi$ has the minimal degree and no presentation
  of the same degree exists. Otherwise, we can:
  
  (a) replace some $\alpha_i$ by $\beta_i$,
  (b) add $\beta_n$,
  
  (c) replace some $\beta_j$ by $\alpha_j$,
  (d) add $\alpha_n$.
  
  The steps (a),(b) decrease the $n$-th
  coordinate by 1, the steps (c),(d) increase it by 1.
  The increasing of the degree is: $n$ for (b),(d),
  $n-2i$ for (a), $2j-n$ for (c).
  If $\chi'=\chi+t\varphi_n$, then to obtain the minimal presentation,
  we apply $t$ times (a) and (b), if $t>0$,
  and we apply $-t$ times (c) and (d), if $t<0$.
  Clearly, there exists a unique sequence of steps giving $\chi$ 
  with the minimal possible degree. Therefore
  the presentation of $\chi$ with the minimal degree is unique. 
  Moreover, from its construction follows that the 
  presentation of $c\chi$ with the minimal degree is just the sum of
  $c$ minimal presentations for $\chi$. This completes the proof.$\Box$

  \vspace{0.3cm}

  {\it Proof of Lemma ~\ref{delta}.}

  Consider the case $G=GL(V)$. Let $\varepsilon_1,\cdots,\varepsilon_n$
  be the weights of $T$ acting on $V$, a basis of the character
  lattice of $T$.
  let $\chi_1,\cdots,\chi_n$ be the dual basis.
  The fundamental weights are:
  $\varphi_i=\varepsilon_1+\cdots+\varepsilon_i, i=1,\cdots, n$.
  Furthermore,  ${\cal C}$ is given by the inequalities
  $\chi_1\geq\chi_2\cdots\geq\chi_n$,
  $\Phi={\rm conv}(\pm\varepsilon_1,\cdots,\pm\varepsilon_n),$
  and $\widetilde{\Delta}$ is the convex hull of
  $$\varepsilon_1,\frac{\varepsilon_1+\varepsilon_2}{2},
  \cdots,\frac{\varepsilon_1+\cdots+\varepsilon_n}{n},
  -\varepsilon_n, \frac{-\varepsilon_n-\varepsilon_{n-1}}{2},\cdots,
  \frac{-\varepsilon_n-\cdots-\varepsilon_1}{n}.$$
  
  For $\chi\in\langle\varepsilon_1,\cdots,\varepsilon_n\rangle_{\bf Q}$,
  set $\alpha_i=\chi_i(\xi)$. First assume
  \begin{equation}
  \label{symplex+} 
  \alpha_1\geq\alpha_2\geq\cdots\geq\alpha_n\geq 0,
  \alpha_1+\cdots+\alpha_n\leq 1.
  \end{equation} 
  Then we can rewrite:      
  $$\xi=(\alpha_1-\alpha_2)\varphi_1+(\alpha_2-\alpha_3)\varphi_2+\cdots           +(\alpha_{n-1}-\alpha_n)\varphi_{n-1}+\alpha_n\varphi_n.$$
  So $\xi$ is a linear combination of 
  {\Large $\frac{\varphi_i}{i}$},$i=1,\cdots,n$
  with non-negative coefficients. Now we sum the coefficients: 
  $$(\alpha_1-\alpha_2)+2(\alpha_2-\alpha_3)+\cdots+
  (n-1)(\alpha_{n-1}-\alpha_n)+n\alpha_n=\alpha_1+\cdots+\alpha_n\leq 1.$$
  Therefore we get: 
  $$\xi\in
  {\rm conv}(0,\varepsilon_1,
  \frac{\varepsilon_1+\varepsilon_2}{2},\cdots,
  \frac{\varepsilon_1+\cdots+\varepsilon_n}{n})
  \subseteq\widetilde{\Delta}.$$
  
  Analogously, assuming 
  \begin{equation}\label{symplex-}
  0\geq\alpha_1\geq\cdots\geq\alpha_n,
  \alpha_1+\cdots+\alpha_n\geq -1,
  \end{equation}
  we obtain 
  $$\xi\in{\rm conv}
  (0,-\varepsilon_n,\frac{-\varepsilon_n-\varepsilon_{n-1}}{2},
  \cdots,\frac{-\varepsilon_n-\cdots-\varepsilon_1}{n})
  \subseteq\widetilde{\Delta}.$$

  Now assume $\xi\in\Phi\cap{\cal C}$. Then $\xi\in\Phi$ implies
  $|\alpha_1|+\cdots +|\alpha_n|\leq 1$. 
  If all the $\alpha_i$ are of the same sign, then either (~\ref{symplex+})
  or (~\ref{symplex-}) holds and we are done. 
  Otherwise for some $q<n$ we have 
  $$\alpha_1\geq\cdots\geq\alpha_q\geq 0\geq \alpha_{q+1}
  \geq\cdots\geq\alpha_n.$$
  Then set:
  $$t=\sum_{i=1}^q\alpha_i\leq 1,
  \xi_+=\frac{\sum_{i=1}^q\alpha_i\varepsilon_i}{t},
  \xi_-=\frac{\sum_{j=q+1}^n\alpha_j\varepsilon_j}{1-t}.$$
  Clearly, (~\ref{symplex+}) holds for $\xi_+$ and
  (~\ref{symplex-}) holds for $\xi_-$. 
  Hence, $\xi_+,\xi_-\in\widetilde{\Delta}$, 
  and $\xi=t\xi_++(1-t)\xi_- \in[\xi_+,\xi_-]\subseteq\widetilde{\Delta}$.

  For $G=Sp(V),O(V)$, we let $\varepsilon_1,\cdots,\varepsilon_r$
  to be basic characters of $T(G)$ and keep the notation of $\chi_i$-s.
  Then the fundamental weights are (see e.g. ~\cite{ov}):

  for $G=Sp(V)$, $\varphi_i=\varepsilon_1+\cdots+\varepsilon_i,$
  for $i=1,\cdots,r$,

  for $G=O(V)$, $n=2r+1$, $\varphi_i=\varepsilon_1+\cdots+\varepsilon_i,$
  for $i=1,\cdots,r-1$,

  $\varphi_r=\frac{1}{2}(\varepsilon_1+\cdots+\varepsilon_r),$

  for $G=O(V)$, $n=2r$, $\varphi_i=\varepsilon_1+\cdots+\varepsilon_i,$
  for $i=1,\cdots,r-2$,

  $\varphi_{r-1}=\frac{1}{2}(\varepsilon_1+\cdots+\varepsilon_r),$
  $\varphi_r=\frac{1}{2}(\varepsilon_1+\cdots
  +\varepsilon_{r-1}-\varepsilon_r).$

  For the cases $G=Sp(V)$, $n=2r$ or $G=SO(V)$, $n=2r+1$, we have:
  ${\cal C}$ is given by the inequalities
  $\chi_1\geq\chi_2\cdots\geq\chi_m\geq 0$,
  $$\Phi={\rm conv}(\pm\varepsilon_1,\cdots,\pm\varepsilon_r),
  \widetilde{\Delta}={\rm conv}
  (0,\varepsilon_1,\frac{\varepsilon_1+\varepsilon_2}{2},
  \cdots,\frac{\varepsilon_1+\cdots+\varepsilon_r}{r}).$$
  Therefore for $\xi\in{\cal C}\cap\Phi$ the assumption
  (~\ref{symplex+}) holds, hence $\xi\in\widetilde{\Delta}$.

  For the case $G=O(V),n=2r$, we have:
  $\Phi={\rm conv}(\pm\varepsilon_1,\cdots,\pm\varepsilon_r)$,
  $$\widetilde{\Delta}={\rm conv}
  (0,\varepsilon_1,\frac{\varepsilon_1+\varepsilon_2}{2},
  \cdots,\frac{\varepsilon_1+\cdots+\varepsilon_r}{r},
  \frac{\varepsilon_1+\cdots+
  \varepsilon_{r-1}-\varepsilon_r}{r}).$$
  If $\xi\in{\cal C}$, then we can write:
  $$\xi=\alpha_1\varepsilon_1+\alpha_2\frac{\varepsilon_1+\varepsilon_2}{2}+
  \cdots+\alpha_r\frac{\varepsilon_1+\cdots+\varepsilon_r}{r}+
  \beta\frac{\varepsilon_1+\cdots+
  \varepsilon_{r-1}-\varepsilon_r}{r},$$
  where $\alpha_1,\cdots,\alpha_r,\beta\geq 0$, $\alpha_r\beta=0$.
  Assume $\xi\in\Phi$. If $\alpha_r=0$, then, taking into account the
  inequality $\chi_1(\xi)+\cdots+\chi_{r-1}(\xi)-\chi_r(\xi)\leq 1$,
  we obtain $\alpha_1+\cdots+\alpha_{r-1}+\beta\leq 1$. Therefore
  $\xi\in\widetilde{\Delta}$. Similarly, we consider the case $\beta=0$.
  This completes the proof of Lemma~\ref{delta}.$\Box$

  \section{Syzygies.}\label{syz}

  Since we found the generators of ${\bf k}[W]^{U(G)}$, a natural
  question is to describe their syzygies. This is a subject
  of the Second Fundamental Theorem of Invariant Theory for
  the linear group $(U(G),V)$. In this section we present
  some results for $G=GL(V)$.
  Of course, syzygies that we present are also
  syzygies for the orthogonal and symplectic cases, if
  the involved generators are.

  Set $U=U(GL)$ and denote by $W_U$ the spectrum of ${\bf k}[W]^U$.
  Moreover, denote by $\pi_{U,W}$ the quotient map
  $\pi_{U,W}:W\to W_U$ corresponding to the inclusion
  ${\bf k}[W]^U\subseteq {\bf k}[W]$.

  For any $p,l\in{\bf N},1\leq p\leq l$, set
  $L={\bf k}^l\oplus\wedge^2{\bf k}^l\oplus\cdots\oplus\wedge^p{\bf
  k}^l$. Let  ${\cal F}_{p,l}$ denote the set of all
  $(q_1,q_2,\cdots,q_p)$ in $L$
  such that for $i=2,\cdots,p$,
  the $i$-vector $q_i$ is decomposable, and
  $Ann(q_{i-1})\subseteq Ann(q_i)$,  where $Ann(q)=\{x\in V\vert q\wedge x=0\}$.
  
  The subset ${\cal F}_{p,l}$ is not closed in $L$.
  In fact, assume $(q_1,\cdots,q_p)\in {\cal F}_{p,l}$ is such that
  $q_2\neq 0$. Then for any $t\in{\bf k}^*$ the collection
  $(tq_1,q_2,\cdots,q_p)$ also belongs to ${\cal F}_{p,l}$.
  But the limit $(0,q_2,\cdots,q_p)$ of such collections does not
  belong to ${\cal F}_{p,l}$. Denote by $\overline{{\cal F}_{p,l}}$
  the Zariski closure of ${\cal F}_{p,l}$,

  Note that the subset ${\cal F}_{p,l}$ is stable under 
  the natural action of the group $GL_l$ on $L$. Therefore
  ${\cal F}_{p,l}$ is acted upon by $GL_l$.

  \begin{theorem}\label{vectors}
  For $W=lV$, set $p={\rm min}\{l,n\},
  q=n-p+1$.
  Consider the rows $u_1,\cdots,u_n$ of the matrix $\overline{V}$
  as the coordinates of some vectors in ${\bf k}^l$.
  Then the map $W\to\overline{{\cal F}_{p,l}}\subseteq L$ taking a tuple of    
vectors
  to the element with coordinates
  $$(u_n,u_{n-1}\wedge u_n,\cdots,u_q\wedge u_{q+1}\wedge\cdots\wedge
  u_n)$$
  is the $GL_L$-equivariant quotient map $\pi_{U,W}$
  and its image is ${\cal F}_{p,l}$.
  \end{theorem}

  {\it Proof.} We only need to prove that the Pl\"ucker
  coordinates of the antisymmetric forms $u_q\wedge\cdots\wedge u_n,
  \cdots,u_{n-1}\wedge u_n, u_n$ generate ${\bf k}[W]^U$. 
  But these are just
  the lower minor determinants of $\overline{V}$ and
  Theorem ~\ref{main} implies Theorem ~\ref{vectors}.
  A different proof of both Theorems for this case
  is as follows. Let the maximal
  unipotent subgroup  $U'\subseteq GL_l$ consist of
  all the strictly upper triangular matrices, in the
  chosen basis of ${\bf k}^l$. It is well known
  (see e.g. ~\cite[3.7]{kr}) that ${\bf k}[W]^{U\times U'}$
  is generated by the left lower minor determinants of $\overline{V}$.
  Therefore the algebra $A$ generated by all the lower
  minor determinants contains ${\bf k}[W]^{U\times U'}$.
  In other words, $A^{U'}=({\bf k}[W]^U)^{U'}$.
  Since $A$ is $GL_l$-stable, we obtain $A={\bf k}[W]^U$.$\Box$

  Thus the syzygies of the set of lower minor determinants
  of $\overline{V}$ are the generators of the ideal in ${\bf k}[L]$
  vanishing on ${\cal F}_{p,l}$. These are the Pl\"ucker
  relations saying that each $q_i$ is decomposable,
  and the incidence relations saying $Ann(q_i)\subseteq Ann(q_j)$
  for any $1\leq i<j\leq p$.

  The syzygies can be written down explicitly. For instance,
  if $i+j\leq p$, then we construct a $(i+j)\times l$
  matrix of the last $i$ rows and
  the last $j$ rows of $\overline{V}$. Clearly,
  any minor determinant of order $i+j$ of such a matrix
  is zero. This is a bilinear  syzygy among the lower
  minor determinants of order $i$ and $j$.

  There is also a $GL_l$-equivariant description of the ideal
  of syzygies, in the form of ~\cite{br2}. 
  For $1\leq i\leq j\leq p$, let $M_{i,j}$ be
  the $GL_l$-stable complementary subspace to
  the highest vector irreducible factor of
  $(\wedge^i{\bf k}^l)^*\otimes(\wedge^j{\bf k}^l)^*\subseteq{\bf k}[L]$,
  if $i<j$, or of $S^2(\wedge^i{\bf k}^l)^*\subseteq{\bf k}[L]$, if $j=i$.
  Let $J$ be the ideal generated by $M_{i,j}$, for
  all $1\leq i\leq j\leq p$.

  \begin{lemma}\label{ideal}
  The ideal in ${\bf k}[L]$ vanishing on ${\cal F}_{p,l}$
  is $J$.
  \end{lemma}

  {\it Proof:} Clearly, we have: $\overline{{\cal F}_{p,l}}=GL_l(L^{U'})$.
  Then by the Theorem of ~\cite[p.382]{br2}, the set of
  zeros of $J$ is $\overline{{\cal F}_{p,l}}$. Moreover, by the same theorem,
  $J$ is radical. This completes the proof.$\Box$

  \begin{corollary}
  All the syzygies are of degree 2.
  \end{corollary}

 Clearly, for arbitrary $l$ and $m$, similar Pl\"ucker and incidence
 relations hold for the left minor determinants of $\overline{V^*}$.

 \begin{theorem}\label{l+m}

 Suppose that $l>0,m>0$ and set $W=lV+mV^*$. Then the ideal of syzygies
 for the generators of ${\bf k}[W]^U$ is generated by
 the Pl\"ucker and the incidence relations for
 the lower minor determinants of $\overline{V}$ and
 for the left minor determinants of $\overline{V^*}$
 if and only if $l+m\leq n$.
 \end{theorem}

{\it Proof:} To prove the "if"  part, it is sufficient
to consider the case $l+m=n$.
Recall that by Theorem ~\ref{generators}, the generators
of ${\bf k}[W]^U$ are
the lower minor determinants of $\overline{V}$,
the left minor determinants of $\overline{V^*}$,
and the elements of the matrix $C=\overline{V^*V}$.
Let $\sum_\alpha a_{\alpha}c^{\alpha}=0$ be a relation
among the generators, where $c^{\alpha}$ is a monomial in
the $C_i^j$-s, $a_{\alpha}$ is a polynomial in the
minor determinants. The assertion of the Theorem
amounts to prove that $a_{\alpha}$ belongs to the
ideal of syzygies, for any $\alpha$.
This will be proven if we check for generic fibers
$F=\pi_{U,lV}^{-1}(\xi),\xi\in{\cal F}_{l,l}$ and
$F^*=\pi_{U,mV^*}^{-1}(\eta),\eta\in{\cal F}_{m,m}$ that
the restrictions of the matrix elements of $C$ to
$F\times F^*$ are algebraically independent. Fix a tuple
of vectors in a generic fiber $F$ such that $\overline{V}$
has the form
$$\left(
\begin{array}{ccc}
\leftarrow & l & \rightarrow \\
 \ast & * & *\\
 \ast & * & *\\
 \hline
 a_1 & 0 & 0\\
 \ast  & \ddots & 0\\
 \ast & * & a_l\\
\end{array}
\right)
$$
with $a_1a_2\cdots a_l\neq 0$ and fix generic elements
of the first $m$ columns of $\overline{V^*}$.
Then, varying the $lm$ elements in the last
$l=n-m$ columns of $\overline{V^*}$, we do not change
the minor determinants and we can obtain any $m\times l$
matrix as $C$.  Thus the "if" part is proven.

The "only if" part. Take $l,m$ such that $1\leq l,m\leq n,l+m>n$
and set $s=l+m-n, r=n-l+1$. Denote by $a_i^j,b_i^j,c_i^j$ the
element in the $i$-th row and the $j$-th column of the matrix
$\overline{V^*},\overline{V}$, and $C$, respectively.
Denote by $\varepsilon^{a\cdots b}$ and $\varepsilon_{a\cdots b}$
the determinant tensors.
In this notation,
$\varepsilon^{i_1\cdots i_m}a_{i_1}^1\cdots a_{i_m}^m$ is the left minor
determinant of order $m$ of $\overline{V^*}$ and
$\varepsilon_{j_1\cdots j_l}b_r^{j_1}\cdots b_n^{j_l}$ is the lower minor
determinant of order $l$ of $\overline{V}$. We claim that
the following relation holds\footnote{This relation with $m=n$ was
indicated to us by E. B. Vinberg.}:
\begin{equation}\label{zacep}
\varepsilon^{i_1\cdots i_m}a_{i_1}^1\cdots a_{i_m}^m
\varepsilon_{j_1\cdots j_l}b_r^{j_1}\cdots b_n^{j_l}=
\end{equation}
$$=\frac{1}{s!}\varepsilon^{i_1\cdots i_m}a_{i_1}^1\cdots a_{i_{r-1}}^{r-1}
\varepsilon_{j_1\cdots j_l}b_{m+1}^{j_{s+1}}\cdots b_n^{j_l}
c_{i_r}^{j_1}\cdots c_{i_m}^{j_s}.$$
To prove this formula, we rewrite the right hand side, using 
$c_i^j=a_i^k b_k^j$:
\begin{equation}\label{yavno}
\frac{1}{s!}\varepsilon^{i_1\cdots i_m}a_{i_1}^1\cdots a_{i_{r-1}}^{r-1}
a_{i_r}^{k_1}\cdots a_{i_m}^{k_s}
\varepsilon_{j_1\cdots j_l}b_{k_1}^{j_1}\cdots b_{k_s}^{j_s}
b_{m+1}^{j_{s+1}}\cdots b_n^{j_l}.
\end{equation}
Let $S(k_1,\cdots, k_s)$ denote the sum of terms in formula
(~\ref{yavno}) with fixed $k_1,\cdots,k_s$.
 Clearly, if $\{k_1,\cdots,k_s\}\neq\{r,r+1,\cdots,m\}$,
 then $S(k_1,\cdots,k_s)=0$.
 Moreover, if $\{k_1,\cdots,k_s\}=\{r,\cdots,m\}$, then
 $S(k_1,\cdots,k_s)$ equals the left hand side of (~\ref{zacep}).

 Therefore, the relation (~\ref{zacep}) holds. Clearly, the right
 hand side is a polynomial in the left minor determinants of order
 $m-s$ of $\overline{V^*}$, the lower minor determinants of
 order $l-s$ of $\overline{V}$, and the matrix elements of $C$.
 It is not hard to check that this relation among the generators
 of ${\bf k}[W]^U$ can not be obtained from  relations of smaller
 degrees.$\Box$

{\bf Remark.}
Theorems ~\ref{vectors},~\ref{l+m} yield an independent
proof of Theorem ~\ref{generators} for the case $l+m\leq n$.
Indeed, we prove in Theorem ~\ref{l+m} that,
in the case $l+m\leq n$, the syzygies among the elements 
of the set ${\cal M}$ are generated by those for
$lV$ and those for $mV^*$. We did not use Theorem ~\ref{generators}
for this.
Hence, by Theorem ~\ref{vectors} (that we also prove independently
of Theorem ~\ref{generators}), ${\rm Spec} A\cong (lV)_U\times(mV^*)_U$.
Since for an action of an algebraic group $H$ on
a normal affine variety $X$, the algebra ${\bf k}[X]^H$ is integrally
closed, ${\rm Spec} A$ is normal.
Furthermore, as we did it for $O(V)$, one can prove
${\bf k}[W]^U_f=A_f$ for all linear $U$-invariants $f$.
Then any  $g\in {\bf k}[W]^U$ gives rise to a rational
function on ${\rm Spec} A$, regular outside the intersection
of the divisors of these linear $U$-invariants.
Since ${\rm Spec} A$ is normal, we get $g\in A$.

\end{document}